\documentclass[a4paper,12pt]{amsart}
\usepackage{amssymb}
\usepackage{dsfont}
\usepackage{ifthen}
\usepackage{graphicx}

\newtheorem{thm}{Theorem}
\newtheorem{cor}{Corollary}
\newtheorem{lem}{Lemma}

\newtheorem{que}{Question}
\newtheorem{rem}{Remark}
\newtheorem{prob}{Problem}
\newtheorem{conj}{Conjecture}
\theoremstyle{definition}
\newtheorem{example}[equation]{Example}

\newcounter {own}
\def\theown {\thesection       .\arabic{own}}

\newenvironment{pf}[1][]{%
 \vskip 3mm
 \noindent
 \ifthenelse{\equal{#1}{}}%
  {{\slshape Proof. }}%
  {{\slshape #1.} }%
 }%
{\qed\bigskip}

\newcounter{alphabet}
\newcounter{tmp}
\newenvironment{Thm}[1][]{\refstepcounter{alphabet}%
\bigskip%
\noindent%
{\bf Theorem \Alph{alphabet}}%
\ifthenelse{\equal{#1}{}}{}{ (#1)}%
{\bf .} \itshape}{\vskip 8pt}

\makeatletter
\newcommand{\Ref}[1]{\@ifundefined{r@#1}{}{\setcounter{tmp}{\ref{#1}}\Alph{tmp}}}
\makeatother

\newcounter{minutes}
\divide\time by 60
\newcounter{hours}
\multiply\time by 60
\addtocounter{minutes}{-\time}

\def\be{\begin{equation}}
\def\ee{\end{equation}}

\newcommand{\bee}{\begin{enumerate}}
\newcommand{\eee}{\end{enumerate}}

\newcommand{\blem}{\begin{lem}}
\newcommand{\elem}{\end{lem}}
\newcommand{\bthm}{\begin{thm}}
\newcommand{\ethm}{\end{thm}}
\newcommand{\bcor}{\begin{cor}}
\newcommand{\ecor}{\end{cor}}
\newcommand{\beg}{\begin{example}}
\newcommand{\eeg}{\end{example}}
\newcommand{\bques}{\begin{que}}
\newcommand{\eques}{\end{que}}
\newcommand{\begs}{\begin{examples}}
\newcommand{\eegs}{\end{examples}}
\newcommand{\bdefe}{\begin{defin}}
\newcommand{\edefe}{\end{defin}}
\newcommand{\bprob}{\begin{prob}}
\newcommand{\eprob}{\end{prob}}
\newcommand{\bei}{\begin{itemize}}
\newcommand{\eei}{\end{itemize}}

\newcommand{\bcon}{\begin{conj}}
\newcommand{\econ}{\end{conj}}
\newcommand{\bcons}{\begin{conjs}}
\newcommand{\econs}{\end{conjs}}
\newcommand{\bprop}{\begin{propo}}
\newcommand{\eprop}{\end{propo}}
\newcommand{\br}{\begin{rem}}
\newcommand{\er}{\end{rem}}
\newcommand{\brs}{\begin{rems}}
\newcommand{\ers}{\end{rems}}
\newcommand{\bo}{\begin{obser}}
\newcommand{\eo}{\end{obser}}
\newcommand{\bos}{\begin{obsers}}
\newcommand{\eos}{\end{obsers}}
\newcommand{\bpf}{\begin{pf}}
\newcommand{\epf}{\end{pf}}
\newcommand{\ba}{\begin{array}}
\newcommand{\ea}{\end{array}}
\newcommand{\beq}{\begin{eqnarray}}
\newcommand{\beqq}{\begin{eqnarray*}}
\newcommand{\eeq}{\end{eqnarray}}
\newcommand{\eeqq}{\end{eqnarray*}}

\def\cc{\setcounter{equation}{0}   
\setcounter{figure}{0}\setcounter{table}{0}}

\begin{document}

\bibliographystyle{amsplain}

%

\title[Elementary considerations for classes of meromorphic univalent functions]{Elementary considerations for classes of meromorphic univalent functions}

\thanks{
File:~\jobname .tex,
          printed: \number\day-\number\month-\number\year,
          \thehours.\ifnum\theminutes<10{0}\fi\theminutes}

\author[S. Ponnusamy]{Saminathan Ponnusamy
}
\address{S. Ponnusamy, Department of Mathematics,
Indian Institute of Technology Madras, Chennai-600 036, India.
}
\email{samy@isichennai.res.in, samy@iitm.ac.in}

\author[K.-J. Wirths]{Karl-Joachim Wirths}
\address{K.-J. Wirths, Institut f\"ur Analysis und Algebra, TU Braunschweig,
38106 Braunschweig, Germany.}
\email{kjwirths@tu-bs.de}


\begin{abstract}
In this article we consider functions $f$ meromorphic in the unit disk. We give an elementary proof
for a condition that is sufficient for the univalence of such functions. This condition simplifies and generalizes
known conditions. We present some typical problems of geometrical function theory and
give elementary solutions in the case of the above functions.
\end{abstract}

\thanks{The first author is currently at the ISI, Chennai Centre.}

\maketitle
\pagestyle{myheadings}
\markboth{S. Ponnusamy and K.-J. Wirths}{Meromorphic univalent functions}
\cc

\section*{Introduction and Main Results}
Let $\mathds{D}\,=\,\{z:\, |z|< 1\}$ and $f$ be meromorphic in  $\mathds{D}$. Much research in the last century has been done concerning
such functions that  are injective or univalent. Many conditions that are sufficient for univalence of analytical or geometrical character
have been found and considered in detail. The present paper is devoted to the following condition of such type that  has been found by
Aksent\'ev in 1958, see \cite{Aks58}, and proved again with another method by Nunokawa and Ozaki in 1972, see \cite{NO}.

\begin{Thm}\label{Th1}
Let $f$ be holomorphic in $\mathds{D}$ and such that $f(0)=f'(0)-1=0$. If for all $z\in\mathds{D}$ the inequality
\begin{equation}\label{f1}
\left|\frac{z}{f(z)}\,-\,z\left(\frac{z}{f(z)}\right)'\,-\,1\right|\,<\,1\end{equation}
is valid, then $f$ is univalent in $\mathds{D}$.
\end{Thm}

The present paper is devoted to the following simple generalization of this theorem for  which we will provide an elementary proof.

\begin{thm}\label{Th1a}
Let $f$ be meromorphic in $\mathds{D}$ and such that $f(0)=f'(0)-1=0$.  If for all $z\in\mathds{D}$ the inequality
\eqref{f1} is valid, then $f$ is univalent in $\mathds{D}$.
\end{thm}
\bpf
The condition (\ref{f1}) implies immediately that $f$ has no zero in $\mathds{D}\setminus\{0\}$. Further, there exists
a function $\Omega,$ holomorphic in  $\mathds{D}$ and satisfying the condition $|\Omega(z)|\,<\, 1$ for $z\in \mathds{D}$ such that
\begin{equation}\label{f2}
\frac{z}{f(z)}\,-\,z\left(\frac{z}{f(z)}\right)'\,-\,1\,=\,\Omega(z).
\end{equation}
The normalization of $f$ implies that $\Omega$ has a zero of multiplicity at least two at the origin.
Therefore we may write $\Omega(z)=z^2\omega(z)$, where $\omega$ is analytic in $\mathds{D}$ such that
$|\omega(z)|\,\leq\, 1$ for $z\in \mathds{D}$. The integration of the differential equation (\ref{f2})
delivers the representation
\begin{equation}\label{f3}
\frac{z}{f(z)}\,=\,1+\,cz\,-\,z\int_0^z\omega(t)\,dt,
\end{equation}
where $c$ is an arbitrary constant. This representation is known since long from the considerations of Theorem \Ref{Th1},
see for example the references in \cite{OPW}. It is clear from the above that $f$ satisfies (\ref{f1}) if and only if it has the
representation (\ref{f3}). The proof of the univalence is a one line proof now. Consider $z_1,z_2 \in \mathds{D},~ z_1\neq z_2$ and
\[\frac{f(z_1)\,-\,f(z_2)}{z_1\,-\,z_2}
\,=\,\frac{1\,-\, \frac{z_1z_2}{z_1-z_2}\int_{z_1}^{z_2}\omega(t)\,dt}{(1+\,cz_1\,-\,z_1\int_0^{z_1}\omega(t)\,dt)(1+\,cz_2\,-\,z_2\int_0^{z_2}\omega(t)\,dt)}.\]
Since
\begin{equation}\label{f4}
\left|\int_{z_1}^{z_2}\omega(t)\,dt\right|\,\leq |z_1\,-\,z_2|,
\end{equation}
we see that $f(z_1)\,-\,f(z_2)\,\neq 0$. This proves the univalence of $f$.
\epf

\br
{\rm
The difference between Theorem \Ref{Th1} and Theorem \ref{Th1a} is the fact that in Theorem \ref{Th1a},
the function $f$ is allowed to have a pole in the unit disk. The proof of Theorem \ref{Th1a} shows that such pole is necessarily
a simple pole.

We want to add the hint that similar meromorphic functions have been considered in \cite{BP}.
}
\er

In what follows, we shall derive some analytic properties of the meromorphic functions satisfying (\ref{f1}). To that end we refine
this condition, namely, we let $\lambda \in (0,1]$ and consider the condition
\begin{equation}\label{f5}
\left|\frac{z}{f(z)}\,-\,z\left(\frac{z}{f(z)}\right)'\,-\,1\right|\,<\,\lambda,~z\in\mathds{D}.
\end{equation}
For more details on the holomorphic functions satisfying such inequalities, we refer again to \cite{OPW} and the references therein.

\begin{thm}\label{th2}
\begin{enumerate}
\item[(a)]  Let $f$ be meromorphic in $\mathds{D}$  and such that $f(0)=f'(0)-1=0$.  If for all $z\in\mathds{D}$ the inequality
\eqref{f1} is valid, then $f$ is continuous in $\overline{\mathds{D}}$ with a possible exception at a zero of \eqref{f3}.
\item[(b)] Let $f$ be meromorphic in $\mathds{D}$ and such that $f(0)=f'(0)-1=0$. If for all $z\in\mathds{D}$  the inequality
\eqref{f5} is valid for some $\lambda \in (0,1)$, then $f$  is univalent in $\overline{\mathds{D}}$.
\end{enumerate}
\end{thm}
\bpf
 For the proof of (a) it is sufficient to recognize that (\ref{f4}) implies that
\[A(z)\,=\,1+\,cz\,-\,z\int_0^z\omega(t)\,dt
\]
is uniformly continuous in $\mathds{D}$. This implies the continuity of this function in $\overline{\mathds{D}}$ and this in turn
implies the continuity of $f$ in  $\overline{\mathds{D}}$ with the exception of the zeros of $A(z)$.

In the proof of (b) we use the fact that
\[A_{\lambda}(z)\,=\,1+\,cz\,-\,\lambda z\int_0^z\omega(t)\,dt
\]
has a continuous extension to  $\overline{\mathds{D}}$ for $\lambda \in [0,1]$ which may be denoted by $A_{\lambda}(z)$, too.
For $z_1,z_2 \in \overline{\mathds{D}}, z_1\neq z_2$, we get as in the proof of Theorem \ref{Th1a} that
\[\left|\frac{z_1\,A_{\lambda}(z_2)\,-\,z_2\,A_{\lambda}(z_1)}{z_1\,-\,z_2}\right|\,\geq\,1\,-\,\lambda | z_1z_2|\,>\,0
\]
for $\lambda \in  [0,1)$. This implies
\[\frac{f(z_1)\,-\,f(z_2)}{z_1\,-\,z_2}\,\neq\,0,
\]
and therefore the assertion of the theorem is proved.
\epf

Another typical problem in geometric function theory is the estimation of the moduli of the Taylor coefficients.
Since in our case, the functions $f$ are holomorphic in a neighbourhood of the origin, we may ask for the upper bounds of the
coefficients of the Taylor expansion in a disk around the origin. For the second coefficient we can give a
satisfactory result using an elementary method from \cite{OPW}. Also, compare with \cite{BP}, too.

\begin{thm}\label{th3}  Let $\lambda \in (0,1]$  and $f$ be meromorphic in $\mathds{D}$ and such that $f(0)=f'(0)-1=0$.
If for all $z\in\mathds{D}$ the inequality \eqref{f5} is valid and
\[f(z)\,=\,z\,+\sum_{n=2}^{\infty}a_nz^n
\]
has no pole in the disk $\{z:\,|z|<p\}, ~p\in (0,1]$,  then the inequality
\begin{equation}\label{f6}
|a_2|\,\leq \,\frac{1\,+\,\lambda p^2}{p}
\end{equation}
is valid. The inequality is sharp and equality is attained if and only if
\begin{equation}\label{f7}
f(z)\,=\,\frac{z}{\left(1\,-\,\frac{e^{i\theta}z}{p}\right)\left(1\,-\,\lambda p e^{i\theta} z\right)},\quad \theta \in [0,2\pi).
\end{equation}
\end{thm}
\bpf From the hypotheses above, we know that $f$ is of the form
\[f(z)\,=\,\frac{z}{1+cz-\lambda z\int_0^z\omega(t)\,dt}.
\]
Obviously, $a_2\,=\,-c$. Now, we assume on the contrary that
\[|a_2|\,>\,\frac{1\,+\,\lambda p^2}{p}.
\]
Then we can define
\[|a_2|\,=\,\frac{1\,+\,\lambda p^2}{pr}
\]
for some $r\in (0,1)$. We consider the function
\begin{equation}\label{f8}
F(z)\,=\,\frac{1}{a_2}\left (1\,-\,\lambda z\int_0^z\omega(t)\,dt \right )
\end{equation}
in the closed disk $\overline{\mathds{D}}_{rp}:=\{z:\,|z|\,\leq\,pr\}$. Then we have the inequality
\[
|F(z)|\,\leq\,\frac{pr(1\,+\,\lambda p^2r^2)}{1\,+\,\lambda p^2}\,<\,pr
\]
which implies that $F$ maps this disk into itself. Secondly, for $z_1, z_2$ in the disk $\overline{\mathds{D}}_{rp}$,
we have
\begin{eqnarray*}
|F(z_1)\,-\,F(z_2)|\,&=&\,\frac{\lambda r p}{1\,+\,\lambda p^2}\left|z_1\int_0^{z_1}\omega(t)\,dt\,-\,(z_2-z_1+z_1)\int_0^{z_2}\omega(t)\,dt\right|\\
\,&\leq&\, \frac{\lambda r p}{1\,+\,\lambda p^2}\left(|z_1|\left|\int_{z_2}^{z_1}\omega(t)\,dt\right|\,+\,|z_2-z_1|\left|\int_0^{z_2}\omega(t)\,dt\right|\right)\\
\,&\leq&\,\frac{\lambda r p(|z_1|+|z_2|)|z_2-z_1|}{1\,+\,\lambda p^2} \\
\,&\leq&\,\frac{2\lambda r^2 p^2 }{1\,+\,\lambda p^2}|z_2-z_1|.
\end{eqnarray*}
Since $2\lambda r^2 p^2/(1\,+\,\lambda p^2)\,<\,1,$
the conditions of Banach's fixed point theorem are fulfilled and thus, the function $F$ has a fixed point in the closed disk $\{z:\,|z|\,\leq\,pr\}$.
Therefore $f$ has a pole in the disk $\overline{\mathds{D}}_{rp}$. This is a contradiction to the conditions of Theorem \ref{th3}.
It is easy to see that the functions (\ref{f7}) satisfy the conditions of the theorem and that for them equality in (\ref{f6}) is attained.
To prove the rest of the assertion, we use a method from \cite{OPW}, where this has been shown in the case $p=1$. Therefore, we have to
consider here only the cases $p<1$.
Let us assume that there exists a function $\omega$ that is not of constant modulus one and is such that the function
\[f(z)\,=\,\frac{z}{1\,-\,\frac{1+\lambda p^2}{p}\,e^{i\theta}z\,-\,\lambda z \int_0^z\,\omega(t)\,dt}
\]
satisfies the conditions of the theorem. According to the above assumption, there exists a number $c \in [0,1)$ such that for $|z|\leq p$ the inequality
\[
\left|\int_0^z\,\omega(t)\,dt\right|\,\leq \,c|z|\]
is valid. We define
\[
r\,=\,\frac{1\,+\,\lambda c p^2}{ 1\,+\,\lambda  p^2}\,\in (0,1).\]
Now, we consider again the function (\ref{f8}) and we get by considerations similar to the above ones that $F$ is a
contractive mapping of the closed disk $\{z:\,|z|\leq pr\}$ into itself. Hence, $F$ has a fixed point in this disk, and therefore
there exists a pole of $f$ in the same disk. This again contradicts the above conditions. Now, it remains to prove that the only
quadratic polynomials of the form
\[
1\,-\,\frac{1+\lambda p^2}{p}\,e^{i\theta}z\,-\,\lambda e^{i\varphi}z^2\]
that have no zero in the disk $\{z:\,|z| < p\}$ are the polynomials
\[1\,-\,\frac{1+\lambda p^2}{p}\,e^{i\theta}z\,+\,\lambda e^{2i\theta}z^2.\]
This is a simple exercise of calculations and thus, we complete the proof.
\epf

\br {\rm
In a neighbourhood of the origin, the functions (\ref{f7}) have the Taylor expansion
\[f(z)\,=\,\sum_{n=1}^{\infty}\frac{1\,-\,\lambda^np^{2n}}{p^{n-1}(1\,-\,\lambda p^2)}\,e^{in\theta}z^n.
\]
We conjecture that this delivers the upper bound
\[
|a_n|\,\leq\,\frac{1\,-\,\lambda^np^{2n}}{p^{n-1}(1\,-\,\lambda p^2)}
\]
for all functions $f$ satisfying the conditions of Theorem \ref{th3}. In the case $\lambda = 1$ and $p = 1$,
this is a special case of the famous Bieberbach conjecture that has been proved by de Branges in \cite{B}.
For $ \lambda = 1$ and arbitrary $p \in [0,1)$ the validity of the Bieberbach conjecture implies that the above conjecture
is true. This has been proved by Jenkins in \cite{J}. For $\lambda \in [0,1)$ and  $p =1$ the truth of the above conjecture for $n = 3$ and $n = 4$
was shown in \cite{OPW1}. However the general conjecture remains open.  We conclude the article by recalling that
the above conjecture is a generalization of the following conjecture from  \cite{OPW1}.
}
\er

\bcon\label{conj1}
Suppose that a holomorphic function $f$ in satisfies the condition \eqref{f5} for some $0<\lambda \leq 1$.
Then $|a_n|\leq \sum_{k=0}^{n-1}\lambda ^k$ for $n\geq 2$.
\econ



\begin{thebibliography}{99}
\bibitem{Aks58} L.~A. Aksent\'ev,
\emph{Sufficient conditions for univalence of regular  functions (Russian)},
\textrm{Izv. Vys\v{s}. U\v{c}ebn. Zaved. Matematika} \textbf{1958} (4)(1958), 3--7.



\bibitem{BP} B. Bhowmik and F. Parveen,
\emph{On a subclass of meromorphic univalent functions},
Complex Var. Elliptic Equ. {\bf 62} (2017), 494--510.

\bibitem{B} L. de Branges,
\emph{A proof of the Bieberbach conjecture},
Acta Math. {\bf 20} (1985), 137--152.

\bibitem{J} J. A. Jenkins,
\emph{On a conjecture of Goodman concerning meromorphic univalent functions},
Michigan Math. J. {\bf 9} (1962), 25--27.

\bibitem{NO} M. Nunokawa and S. Ozaki,
\emph{The Schwarzian derivative and univalent functions},
Proc. Amer. Math. Soc. {\bf 33} (1972), 392--394.

\bibitem{OPW} M. Obradovi\'{c}, S. Ponnusamy, and K.-J. Wirths,
\emph{Geometric studies on the class $\mathcal{U}(\lambda)$},
Bull. Malaysian Math. Sci. Soc.  {\bf 39} (2016), 1259--1284.

\bibitem{OPW1} M. Obradovi\'{c}, S. Ponnusamy, and K.-J. Wirths,
\emph{Logarithmic coefficients and a coefficient conjecture for univalent functions},
Monatshefte Math. (2017). \\
Online DOI 10.1007/s00605-017-1024-3


\end{thebibliography}
\end{document}